\documentclass[12pt,a4paper]{article}
\usepackage[table]{xcolor}
\usepackage{diagbox}
\usepackage{amsmath}
\usepackage{cancel}
\usepackage{graphicx}
\usepackage{subcaption}
\usepackage{amsmath,amssymb,amsthm,graphics,amsfonts}
\usepackage{multirow}
\usepackage{csquotes}

\nonstopmode\numberwithin{equation}{section}

\newtheorem{Lemma}{Lemma}[section]
\newtheorem{Theorem}[Lemma]{Theorem}
\newtheorem{Conjecture}[Lemma]{Conjecture}
\newtheorem{theorem}{Theorem}

\renewcommand{\qed}{\rule{6pt}{6pt}}
\newenvironment{Proof}{\noindent{\bf Proof:}}{\qed\bigskip}
 
 \usepackage{lipsum}
\usepackage{color}
\newcommand{\TB}{\textcolor{blue}}

 \newcommand{\bea}{\begin{eqnarray}}
 \newcommand{\eea}{\end{eqnarray}}

\usepackage{graphicx,latexsym}
\usepackage{lscape}

\begin{document}
\title{\LARGE Monotonicity of limit wave speed of periodic traveling wave solutions via Abelian integral }
\author{ \Large{
Krishna Patra, 
Ch. Srinivasa Rao}\\
Department of Mathematics, Indian Institute of Technology Madras\\
Chennai 600036, India}
\date{}
\maketitle
\begin{center}{\bf{Abstract}}
\end{center}
In this article, we investigate monotonicity of limit wave speed of periodic traveling wave solutions for a perturbed generalized KdV
equation via Abelian integral. We have answered an open problem outlined by Yan et al. (2014) and the conjecture proposed by Ouyang et al. (2022). Geometric singular perturbation theory allows for the reduction of a three-dimensional dynamical system to a near-Hamiltonian planar system. Furthermore,  utilizing the monotonic behavior of the ratio of Abelian integrals, we develop a method to show the existence of at most one isolated periodic traveling wave which is much simpler proof than that in Yan et al.(2014). Finally, we present numerical simulations that perfectly match the theoretical outcomes.
\section{Introduction}
Shallow water wave hydrodynamics play a crucial role in engineering applications, particularly in the field of ocean engineering, which involves building marine structures and managing coastal zones. Due to its importance, many scholars have developed mathematical models to analyze the behavior of waves. Korteweg–de Vries (KdV) equations are notable models for studying shallow water waves. They capture the unidirectional movement of these waves, whereas Boussinesq-type equations \cite{ref16} are used to represent bi-directional wave propagation. 

The Korteweg–de Vries class of equations holds great importance as a nonlinear evolution model, finding applications in diverse fields such as plasma physics \cite{ref17}, geophysical fluid dynamics \cite{ref19}, and blood circulation in vessels \cite{ref20}.

Recently, numerous researchers have focused on studying shallow water wave models and their traveling wave solutions in the field of fluid mechanics. KdV equation serves a crucial role in the modeling of shallow water dynamics. The simplest Korteweg-de Vries (KdV) equation is 
\begin{equation} \label{eqn1.1}
 U_t+U U_x+U_{xxx}=0.   
\end{equation}
It describes the one-dimensional propagation of small-amplitude, weakly dispersive waves. The nonlinear term $U U_x$ results in wave steepening, while the linear dispersion term $U_{xxx}$ causes the wave to spread. Considering real-world scenarios, it is essential to incorporate some weak terms into the equation to address external uncertainties like weakly dissipative effects \cite{ref21}. For this reason, the perturbed equation should be taken into consideration. The perturbed generalized KdV equation for modeling shallow water on an inclined, thinner, and relatively long layer was proposed by Derks and van Gils \cite{ref22}, and later by Ogawa \cite{ref23} given by
\begin{equation} \label{kpr3.1}
U_t+U U_x+U_{xxx}+\epsilon(U_{xx}+U_{xxxx})=0,
\end{equation}
 where $\epsilon>0$ is sufficiently small. For $\epsilon=0,$ the equation (\ref{kpr3.1}) reduces to the KdV equation (\ref{eqn1.1}).
 
 In this paper, we study the following perturbed generalized KdV equations
 \begin{equation} \label{kpr3.2}
U_t+U^n U_x+U_{xxx}+\epsilon(U_{xx}+U_{xxxx})=0,
\end{equation}
Yan et al.\cite{ref1} showed in 2014 that solitary and periodic wave solutions of the equation (\ref{kpr3.2}) persist when subjected to sufficiently small perturbation parameters. For $n=1$, Ogawa \cite{ref23} established that the limit wave speed $c_{0}(h)$ behaves as a smooth, decreasing function.  For $n=2,3$, the authors in \cite{ref24} confirmed that $c_{0}(h)$ maintains this smooth, decreasing nature. For $n=4,$  Chen et. al.\cite{ref2} further showed that $c_{0}(h)$ exhibits this same behavior. For $n \geq 5,$ the problem remains open. Motivated by the work of \cite{ref23}, \cite{ref2} and \cite{ref24}, we prove that $c_0(h)$ is a monotone function for all positive integers $n.$ 

Recently, numerous researchers have focused on reaction–convection–diffusion equations due to their significant applications in various scientific and engineering fields. These equations describe intriguing phenomena such as fluid flow, heat transfer, and population dynamics. Several authors have recently proven the existence of periodic traveling wave solutions for different types of reaction–convection–diffusion equations(see, for instance, \cite{ref26}, \cite{ref4}, \cite{ref29}, \cite{ref30},  \cite{ref12}, \cite{ref27}, \cite{ref28},  \cite{ref25}). In a recent study, Fan and Wei \cite{ref33} analyzed traveling waves within a quintic BBM equation under the influence of distributed delay and weak backward diffusion, where they proved the uniqueness of periodic waves by studying the monotonicity properties of Abelian integrals.

The remainder of the paper is structured as follows.
Section 2 covers  an overview of preliminary discussions and utilizes geometric singular perturbation theory and regular perturbation analysis to examine the existence of periodic wave solutions for equation (\ref{kpr3.2}). In Section 3, we first analyze the existence and uniqueness of the periodic wave solution for (\ref{kpr3.2}) for sufficiently small $\epsilon$ with $0<\epsilon\ll 1$ by examining the monotonicity of the ratio of Abelian integrals, after which we prove the monotonicity of the limit wave speed. In Section 4, numerical simulations are provided to demonstrate the existence of a limit cycle and examine the monotonic properties of the limit wave speed for periodic waves. Finally, the paper concludes with Section 5. 

\section{Perturbation Analysis}
Let us take a traveling wave transformation
\begin{equation}\label{kpr2.2}
U(x,t) =U(\xi),\hspace{1cm} \xi = x-ct.
\end{equation}
Here $c$ indicates the wave speed. 
Now substituting (\ref{kpr2.2}) into (\ref{kpr3.2}), we get
\begin{equation}\label{kpr2.3}
 -c \frac{dU}{d\xi}+U^n \frac{dU}{d\xi}+\frac{d^3U}{d\xi^3} +\epsilon\left(\frac{d^2U}{d\xi^2}+\frac{d^4U}{d\xi^4}\right)=0.
\end{equation}
After integrating this equation and choosing the constant of integration as zero, we have
\begin{equation}\label{kpr3.3}
 -c U+\frac{1}{n+1}U^{n+1}+\frac{d^2U}{d\xi^2} +\epsilon\left(\frac{dU}{d\xi}+\frac{d^3U}{d\xi^3}\right)=0.
\end{equation}
Furthermore, applying the scale transformations 
$U=\sqrt[n]{c}u $ and $\xi=\frac{\eta}{\sqrt{c}}$
to the equation (\ref{kpr3.3}), we get
\begin{equation}\label{kpr3.4}
 -u+\frac{1}{n+1}u^{n+1}+\frac{d^2u}{d\eta^2} +\epsilon\left(\frac{1}{\sqrt{c}}\frac{du}{d\eta}+\sqrt{c}\frac{d^3u}{d\eta^3}\right)=0.
\end{equation}
The periodic orbits of system (\ref{kpr3.4}) corresponds to the periodic traveling wave solutions of equation (\ref{kpr3.2}). Therefore, we focus on analyzing the corresponding periodic orbits in equation (\ref{kpr3.4}) in order to explore the periodic waves in equation (\ref{kpr3.2}). Due to the difficulty of directly analyzing equation (\ref{kpr3.4}), we apply geometric singular perturbation theory \cite{ref6} to simplify it onto a slow manifold, where it becomes a regular perturbed problem.


Now equation (\ref{kpr3.4}) can be written as \\
\begin{equation}\label{eqn2.6}
\begin{cases}
     \frac{du}{d \eta}=y,\\ \\ \frac{dy}{d \eta}=z, \\ \\
      \epsilon \sqrt{c} \frac{d z}{d \eta}=u-\frac{1}{n+1} u^{n+1} -z -\frac{\epsilon}{\sqrt{c}}y.
    \end{cases}     
\end{equation}
In \cite{ref1}, Lemma 1, which is based on the framework of geometric singular perturbation theory as introduced by Fenichel \cite{ref6}, demonstrates the existence of the slow manifold  $M_{\epsilon}.$ The detailed derivation and existence of the manifold $M_{\epsilon}$ are thoroughly addressed in \cite{ref1}.

 Now, the dynamics of (\ref{eqn2.6}) on the manifold $M_{\epsilon}$ and omitting $\mathcal{O}(\epsilon^2)$ is given as
\begin{equation}{\label {a01}}
\frac{du}{d\eta}=y,\hspace{.5cm}
 \frac{dy}{d\eta}= u-\frac{1}{n+1}u^{n+1}+\epsilon \sqrt{c}(u^n y-(1+\frac{1}{c})y).
\end{equation}
 Now this is a regular perturbation problem. Now the corresponding unperturbed system [$\epsilon=0$] is 
\begin{equation}{\label {a02}}
\frac{du}{d\eta}=y,\hspace{.5cm}
 \frac{dy}{d\eta}= u-\frac{1}{n+1}u^{n+1}.
\end{equation}
 Now the system (\ref{a02}) has two fixed points $(0,0)$ and $(\sqrt[n]{n+1},0)$ when $n$ is odd, and three fixed points $(0,0)$, $(-\sqrt[n]{n+1},0)$ and $(\sqrt[n]{n+1},0)$ when $n$ is even. The system (\ref{a02}) is a Hamiltonian system \cite{ref32} with Hamiltonian function
 \begin{equation} \label{eqn2.10}
     H(u,y)=\frac{y^2}{2}-\frac{u^2}{2}+\frac{u^{n+2}}{(n+1)(n+2)}:=\frac{y^2}{2}+\Phi(u).
 \end{equation}
  At the fixed points, the Hamiltonian function $H(u,y)$ has values
 \begin{equation}
H(\sqrt[n]{n+1},0)=H(-\sqrt[n]{n+1},0)=\frac{-n(n+1)^{\frac{2}{n}}}{2(n+2)},\hspace{.5cm}H(0,0)=0.
\end{equation}
The theorem below plays a significant role in identifying a fixed point's characteristics as a center.

\begin{theorem}(Hale \cite{ref36}){\label{kp.r2}}
Suppose $\textbf{f}:\mathbb{R}^2\to \mathbb{R}^2$ is a continuous function and for any 
$\textbf{x}_0 \in \mathbb{R}^2$ the system 
 \begin{equation} 
 \dot{\textbf{x}}=\textbf{f}(\textbf{x}),
\label{a1}
\end{equation} 
has a unique solution $\textbf{x}=\textbf{x}(t,\textbf{x}_0)$ with $\textbf{x}(0,\textbf{x}_0)=\textbf{x}_0.$ Assume that the system (\ref{a1}) is conservative and $\textbf{a}\in \mathbb{R}^2$ is an equilibrium point of the system (\ref{a1}). \\
If $E(\textbf{x})$ is an integral of (\ref{a1}) in a bounded open neighborhood $U\subset \mathbb{R}^2$ of the equilibrium point  $\textbf{a}$ such that $E(\textbf{x})> E(\textbf{a})$ for  $\textbf{x} \in U\setminus{\{\textbf{a}\}},$ then the equilibrium point  $\textbf{a}$ is a center. 
\end{theorem}
We start by performing a linear analysis around the fixed points. Now when $n$ is even, the system (\ref{a02}) is symmetric about $u-$ axis and $y-$ axis. So nature of the fixed point $(-\sqrt[n]{n+1},0)$ will be same as nature of the fixed point $(\sqrt[n]{n+1},0).$ Now linear analysis suggests that $(0,0)$ is a saddle for both odd and even $n.$ Also $(\pm \sqrt[n]{n+1},0)$ are center or focus for even $n$, and $(\sqrt[n]{n+1},0)$ is a center or focus for odd $n.$ We will check nature of the point $(\sqrt[n]{n+1},0)$ for both cases. 

 Note that
\begin{eqnarray}
&&H_{u}(u,y)=-u+\frac{u^{n+1}}{n+1},\,\,H_y(u,y)=y, \nonumber\\
&&H_{u y}(u,y) =H_{y u}(u,y) =0, \,\,H_{u u}(u,y)=-1+u^{n}\,\, \TB{,~ } H_{ y y}(u,y)=1, \nonumber\\
&&\tilde{H}(u,y):= H_{u u}(u,y) H_{ y y}(u,y)-(H_{u y}(u,y))^2=u^{n}-1,\nonumber 
\end{eqnarray}
where $H$ is defined in (\ref{eqn2.10}). Now at $(\sqrt[n]{n+1},0),$ $\tilde{H}(\sqrt[n]{n+1},0)= n>0$ for both odd and even $n.$ This clearly indicates that $(\sqrt[n]{n+1},0)$ is an isolated minimum of the function $H(u,y).$ Therefore, by Theorem \ref{kp.r2}, $(\sqrt[n]{n+1},0)$ is a center for the system (\ref{a02}) for both odd and even $n.$ Additionally, $(-\sqrt[n]{n+1},0)$ is a center for even $n$ as nature of the fixed point $(-\sqrt[n]{n+1},0)$ is same as nature of the fixed point $(\sqrt[n]{n+1},0).$ Nature of all fixed points of the unperturbed system (\ref{a02}) are tabulated in Table \ref{tab1}.\\
 
\textbf{$h-$intervals for periodic annuli:}\\ \\
Here, we identify the interval $(p_1,p_2)$ such that the family  $\{\Gamma_h:H(u,y)=h,h \in (p_1,p_2)\}$ of trajectories of (\ref{a02}) forms a periodic annulus  around the center $(\sqrt[n]{n+1},0).$
Note that $H(\sqrt[n]{n+1},0)=-\frac{n(n+1)^{\frac{2}{n}}}{2(n+2)}$ and $H(\sqrt[n]{n+1},y)=\frac{y^2}{2}-\frac{n(n+1)^{\frac{2}{n}}}{2(n+2)}.$ Assuming that the line $u=\sqrt[n]{n+1}$ intersects the trajectory $H(u,y)=h$ of (\ref{a02}). Then, 
solving for $y$ in the equation $H(\sqrt[n]{n+1},y)=h$ yields
$$
y =\pm \sqrt{2\left\{h+\frac{n(n+1)^{\frac{2}{n}}}{2(n+2)}\right\}}. 
$$
This suggests that for $h >  -\frac{n(n+1)^{\frac{2}{n}}}{2(n+2)},$ the trajectories $H(u,y)=h$ of (\ref{a02}) intersect the line $u=\sqrt[n]{n+1}$ at two distinct points. The trajectory $H(u,y) =-\frac{n(n+1)^{\frac{2}{n}}}{2(n+2)} $ corresponds to fixed point $(\sqrt[n]{n+1},0).$ Note that $(0,0)$ is a saddle point, with $H(0,0)=0$. Clearly, the trajectory $H(u,y)=0$ exhibits symmetry about the $u-$axis and intersects $u-$axis at $(0,0)$ and $\left(\left(\frac{(n+1)(n+2)}{2}\right)^{\frac{1}{n}},0\right).$ Moreover, the trajectory $H(u,y)=0$ forms a homoclinic orbit that connects to $(0,0).$ 
Therefore, for $h >  -\frac{n(n+1)^{\frac{2}{n}}}{2(n+2)}$, the trajectory $H(u,y)=h$ expands as $h$ increases,  and the trajectory $H(u,y)=0 $ corresponds to the homoclinic orbit connecting $(0,0)$ and passing through $\left(\left(\frac{(n+1)(n+2)}{2}\right)^{\frac{1}{n}},0\right).$ Hence, the trajectories $\Gamma_h: H(u,y)=h,\hspace{.1cm} h\in \left(-\frac{n(n+1)^{\frac{2}{n}}}{2(n+2)},0\right)$ form a periodic annulus around $(\sqrt[n]{n+1},0).$ So

$$p_1=-\frac{n(n+1)^{\frac{2}{n}}}{2(n+2)},~p_2=0. $$

The periodic orbits of system (\ref{a01}) correspond to periodic traveling wave solutions of equation (\ref{kpr3.2}). Under perturbation, most of the periodic orbits in (\ref{a02}) are broken, with only a few remaining as limit cycles in system (\ref{a01}).

\begin{figure}[htbp] 
  \centering
      \includegraphics[width=1.08\linewidth]{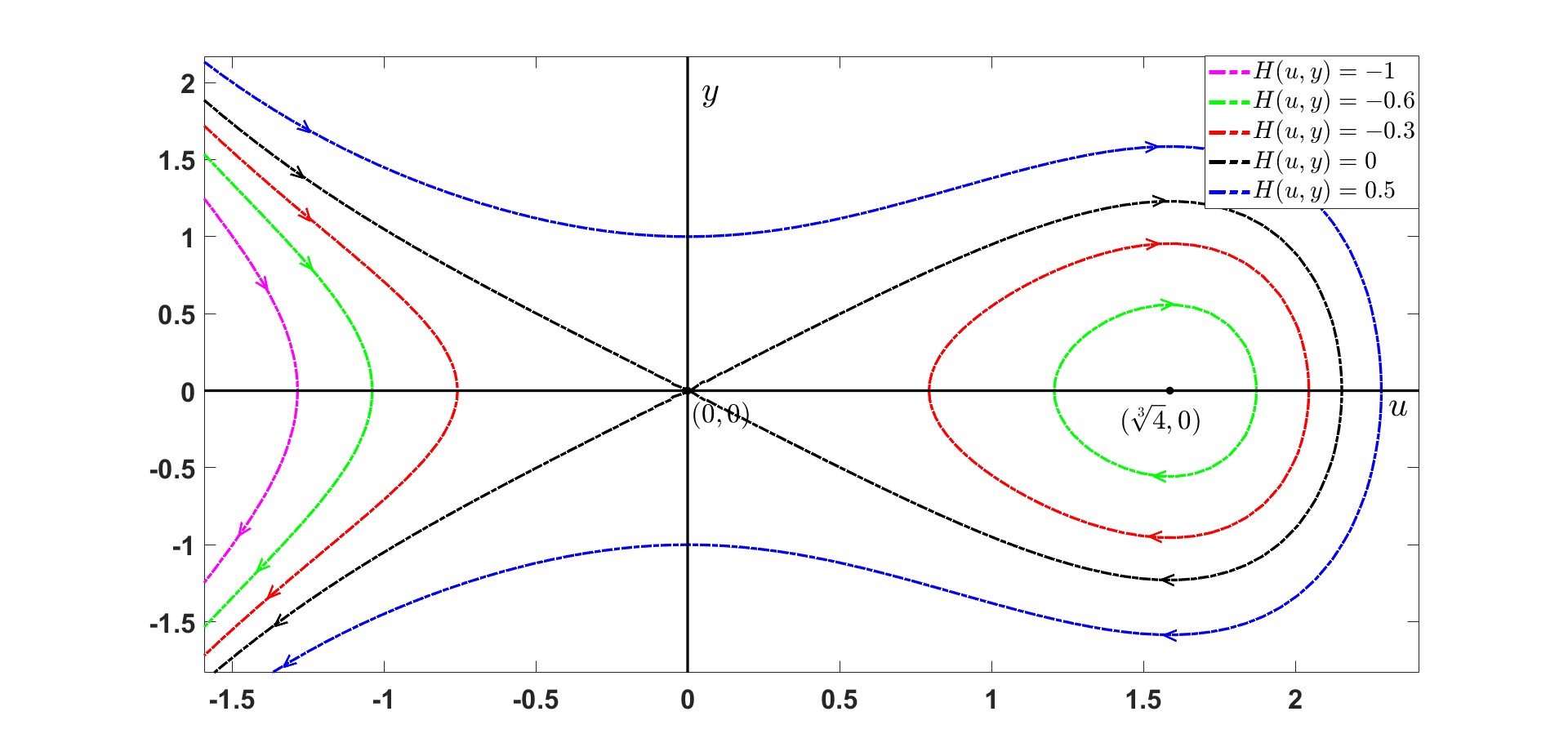}
  \caption{Phase portrait of the system (\ref{a02}) when $n$ is odd($n=3).$}
  \label{kpr001}
\end{figure}

\begin{figure}[htbp] 
  \centering
      \includegraphics[width=1.08\linewidth]{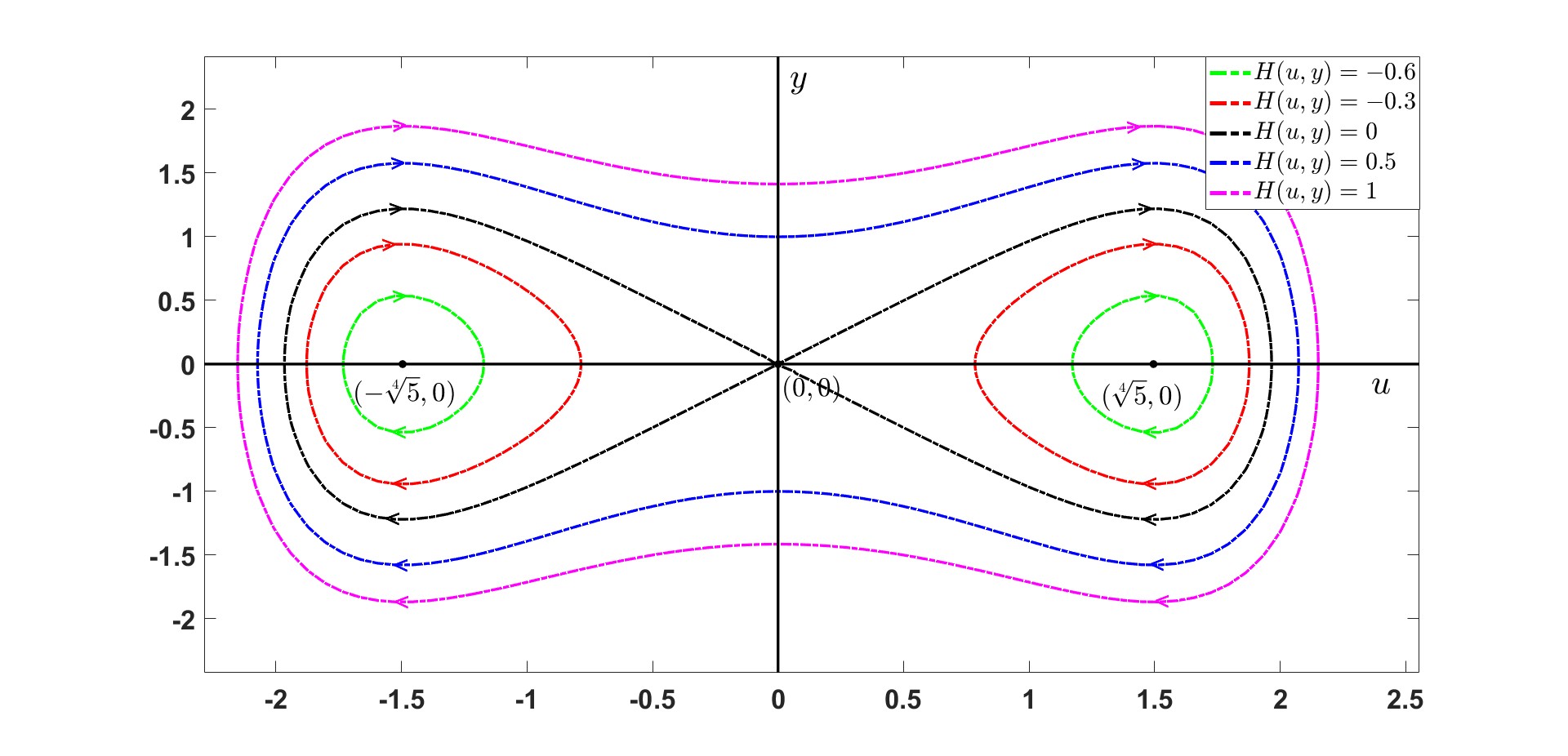}
  \caption{Phase portrait of the system (\ref{a02}) when $n$ is even($n=4).$}
  \label{kpr002}
\end{figure}

Abelian integrals play a significant role in identifying the number of limit cycles present in a perturbed system.

According to the Poincar\'{e}-Pontryagin theorem (as shown in Poincar\'{e} \cite{ref40}  and Pontryagin \cite{ref41}) and Theorem (2.4) in \cite{ref8} indicate that if $\bar{A}(h)$ is not identically zero, the number of its zeros provides an upper estimate for the number of limit cycles in the system (\ref{a01}) that can bifurcate from the periodic annulus.

Now the Abelian integral corresponding to the above perturbed system (\ref{a01}) is
\begin{eqnarray}
 A(h) &=& \oint_{\Gamma_ h}\,\sqrt{c}\left(u^n y-\left(1+\frac{1}{c}\right)y\right)\, du\nonumber\\
 &=&-\sqrt{c}\left(1+\frac{1}{c}\right)\left(\oint_{\Gamma_ h}\,y\, du\right) +\sqrt{c}\left(\oint_{\Gamma_ h}\,u^n y\, du\right)\nonumber\\
&:=& -\sqrt{c}\left(1+\frac{1}{c}\right) A_0(h)+\sqrt{c} A_n(h)\nonumber\\
&=& -\sqrt{c}A_0(h) \left\{\left(1+\frac{1}{c}\right) -\frac{A_n(h)}{A_0(h)} \right\}. \label{eqn2.11}
\end{eqnarray}
See Zoladek and Murilo \cite{ref31}, Patra and Rao \cite{ref4} for more detail explanations about Abelian Integral. Now $ A_0(h)=\oint_{\Gamma_ h}\,y\, du=\iint_{\tilde D}dudy=$Area of $\tilde D>0$ by Green's Theorem in the plane, where $\tilde D$ is the region enclosed by ${\Gamma_h}.$ So we can consider the ratio $\frac{A_n(h)}{A_0(h)}.$ Now if we can show that the ratio $\frac{A_n(h)}{A_0(h)}$ is monotone then we can say that there is at most one periodic wave solution of (\ref{kpr3.2}). Also this information helps us to get desired results. Now the below theorem is helpful to prove monotonicity of the ratio $F_{n}(h):=\frac{A_n(h)}{A_0(h)}$. 
\begin{theorem}{(Liu et al. \cite{ref9})}{\label {thmc}}
Consider the system of differential equations
\begin{equation}{\label {eqn1.13}}
\frac{du}{d\eta}=y,\hspace{.5cm}
 \frac{dy}{d\eta}=\epsilon(\alpha_0 +\alpha_n \tilde g_n(u))y- \Phi'(u), 
\end{equation}
where $\epsilon>0$ is sufficiently small, $\tilde g_n(u)$ is a polynomial, and $\alpha_0,\,\alpha_m$ are real numbers. \\
Suppose that: 
\begin{itemize}
\item[(a)] The Hamiltonian function $H(u,y)$ corresponding to the unperturbed system of $(\ref{eqn1.13})$ has the form $H(u,y) := \Psi(y)+\Phi(u)$ and satisfies $\Phi'(u)(u-\hat a)>0\,\, (\,\, { or} < 0\,\,)$ { for} $u\in(\alpha,B)\setminus \{\hat a\}$ where 
$(\hat a,0)$ is a center of the unperturbed system of $(\ref{eqn1.13})$ and $B$ can be determined by the relation $\Phi(\alpha)=\Phi(B).$\\
(Note that for the condition $(a)$, there exists an involution $\delta$ defined on $(\alpha,B)).$
\item[(b)]
$\hat{A}_{n}(h):=\frac{A_n(h)}{A_0(h)}=\frac{\oint_{\Gamma_h}\,\tilde g_n(u)y\, du}{\oint_{\Gamma_h}\,y\, du}$ and $T_n(u):=(n+1)\frac{\int_{\delta(u)}^{u} \tilde g_n(t) \,dt }{\int_{\delta(u)}^{u} \,dt }$, where $\Gamma_h$ is a closed trajectory of the unperturbed system of $(\ref{eqn1.13})$ and  $\Gamma_h: H(u,y)=h$ 
 for some $h \in (p_1,p_2)$, where $(p_1,p_2)$ is the maximal interval so that the family $\{\Gamma_h:H(u,y)=h,h \in (p_1,p_2)\}$ forms a periodic annulus, and $\delta$ is an involution  defined on $(\alpha,B)).$
\end{itemize}
Then 
 ${T_n}'(u)>0\,\, ( \, \, resp. <0)$ in the interval  $(\hat a,B)$ implies that $\hat{A}'_{n}(h) >0\,\, (resp. <0)$, for all $h \in (p_1, p_2).$
\end{theorem} 

\noindent \textbf{Remark:} As noted in Lemma 2.2 and Remark 2.1 of \cite{ref50}, the sign of $\bar{\xi}'(x)$ in the conclusion of Theorem 2.1 in \cite{ref9} was opposite. The correct version should read: $\bar{\xi}'(x)>0$ (resp.$<0$) in $(0,\nu)$ implies $u'(h)>0$ (resp.$<0$). \\

Also below theorem is helpful to prove existence/non-existence of limit cycle.

\begin{theorem}(Christopher and Li \cite{ref8})\label{lemma3.1}
Assume that $A(h) \not\equiv0$ for $h\in(p_1,p_2)$. Then the following statements hold.\\
\rm{(i)} If a limit cycle of the system $X_{H,\epsilon}$ bifurcates from $\Gamma_{\hat h}$, then $A(\hat h)=0.$\\
\rm{(ii)} If there exists  $\tilde h\in(p_1,p_2)$ such that $A(\tilde h)=0$ and $A'(\tilde h)\neq0$. Then the system $X_{H,\epsilon}$ has a unique limit cycle that bifurcates from $\Gamma_{\tilde h}.$
\end{theorem}

\begin{table}[!h]
\begin{center}
    \begin{tabular}{ |p{1cm} |p{2cm} |p{2cm} | p{2cm} |}
    \hline 
    $~~n$ & ~~$(0,0)$& $(\sqrt[n]{n+1},0)$&$(-\sqrt[n]{n+1},0)$ \\ \hline 
    ~~Even & ~~Saddle&~~\TB{Center}&~~\TB{Center} \\ \hline 
    ~~Odd & ~~Saddle&~~\TB{Center}&~~~~~~-\\ \hline 
    \end{tabular}
  
\end{center}
  \caption{Classification of the fixed points of the unperturbed system (\ref{a02}).}
\label{tab1}
\end{table}

\section{Main results}
\begin{figure}[htbp] 
  \centering
      \includegraphics[width=1.08\linewidth]{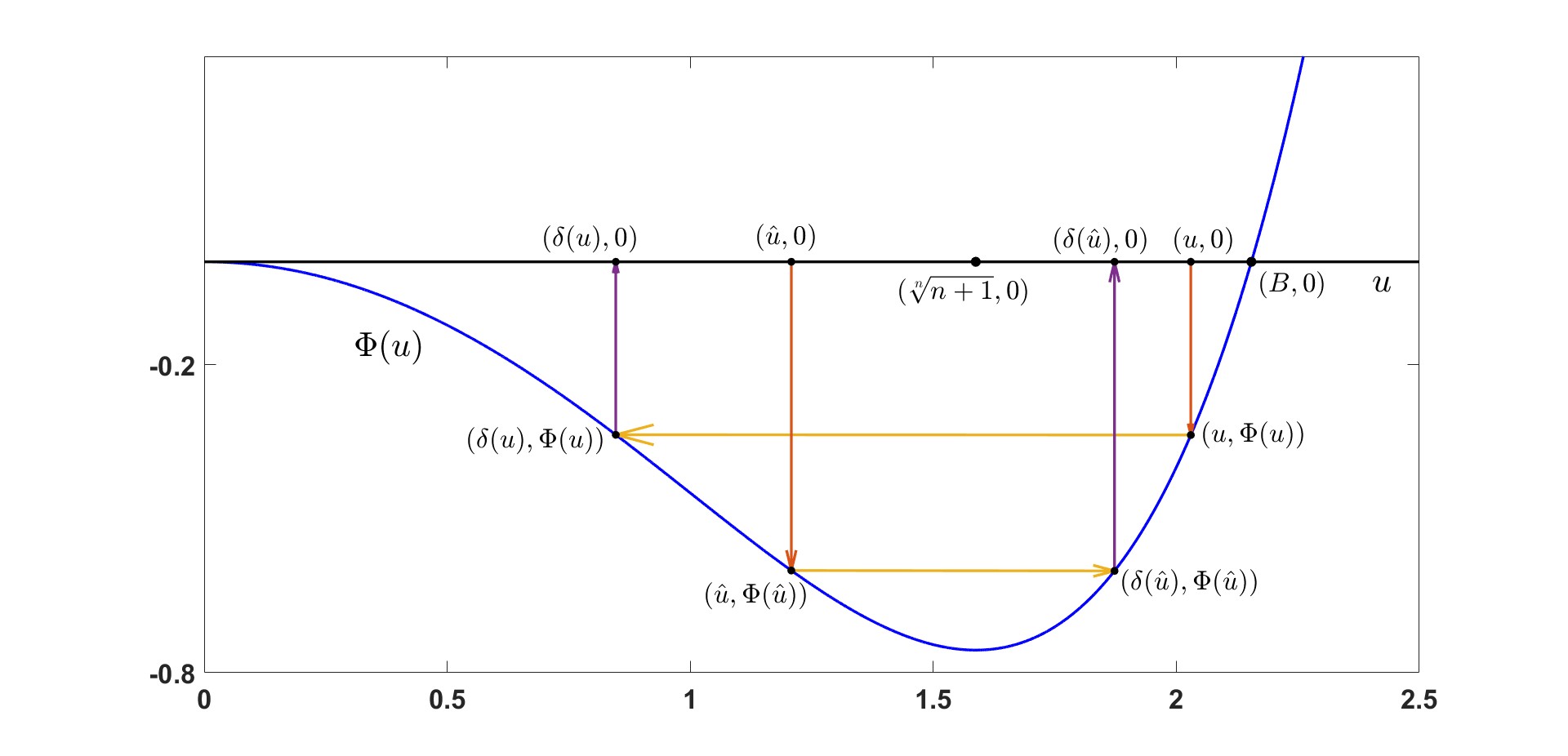}
  \caption{Plot of image of $\Phi$ and involution. }
  \label{kpr003}
\end{figure}

We have shown that $(\sqrt[n]{n+1},0)$ is a center of the system (\ref{a02}) for both even and odd $n$. We can verify that $\Phi'(u)(u-\sqrt[n]{n+1})>0$ for $u\in (0,B)\setminus \{\sqrt[n]{n+1}\},$ where $B$ can be evaluated by the relations $\Phi(B)=\Phi(0)$ and $B>\sqrt[n]{n+1}.$ It is easy to note that $\Phi |_{(0,\sqrt[n]{n+1})}: (0,\sqrt[n]{n+1}) \rightarrow (p_1, \,p_2)$ is strictly decreasing, and $\Phi |_{(\sqrt[n]{n+1},B)}: (\sqrt[n]{n+1},B) \rightarrow (p_1, \,p_2)$ is strictly increasing. Suppose that  $\mu$ and $\nu$ are the inverse functions of $\Phi |_{(0,\sqrt[n]{n+1})}$ and $\Phi |_{(\sqrt[n]{n+1},B)},$ respectively. Then $\mu: (p_1, p_2) \rightarrow (0,\sqrt[n]{n+1})$,  
$\nu: (p_1,p_2) \rightarrow (\sqrt[n]{n+1}, B)$ and 
$$\Phi(\mu(\omega)) \equiv \omega \equiv \Phi(\nu(\omega)) , \,\, \omega\in (p_{1},p_{2}).$$
Here $0<\mu(\omega)<\sqrt[n]{n+1}<\nu(\omega)<B.$ We can easily obtain that there exists an involution $\delta$ defined on $(0,B):$ 
\begin{equation}
    \delta(u)= 
\begin{cases}
    (\nu \circ \Phi)(u) ,& \text{when } u\in (0,\sqrt[n]{n+1});\\
    (\mu \circ \Phi)(u),              & \text{when } u\in  (\sqrt[n]{n+1},B);\\
    \sqrt[n]{n+1},&\text{when } u=\sqrt[n]{n+1},
\end{cases}
\label{inv}
\end{equation}
such that  $\Phi(u)=\Phi(\delta(u)).$ See figure \ref{kpr003} for geometrical observation. Let $v:=\delta(u).$ Clearly if $u\in (0,\sqrt[n]{n+1})$, then $v \in (\sqrt[n]{n+1},B)$ and if $u \in (\sqrt[n]{n+1},B)$ then $v\in (0,\sqrt[n]{n+1}).$ See figure \ref{kpr003}. Here we will take $u \in (\sqrt[n]{n+1},B)$ as similar analysis can be done for other case. So we have 
\begin{equation} \label{eqn3.01}
    0<v<\sqrt[n]{n+1}<u<B.
\end{equation}

Now 
\begin{eqnarray}
&&\Phi(u)=\Phi(v) \nonumber \\
&\implies& 
-\frac{u^2}{2}+\frac{u^{n+2}}{(n+1)(n+2)}=-\frac{v^2}{2}+\frac{v^{n+2}}{(n+1)(n+2)} \nonumber\\
&\implies&u^2-v^2= \frac{2}{(n+1)(n+2)}(u^{n+2}-v^{n+2}). \label{eqn4.2}
\end{eqnarray}
Next Lemma is helpful to prove the conjecture.
\begin{Lemma} \label{kp1a}
 Suppose that $n$ is a positive integer. Then  
\begin{equation} 
(u^2-v^2)\sum_{i=0}^{[\frac{n-1}{2}]}(u^{n-2i}-v^{n-2i})^2(uv)^{2i}=(u^{2(n+1)}-v^{2(n+1)})-(n+1)u^n v^n(u^2-v^2).
\end{equation}
 \end{Lemma}

\begin{Proof}
  \textbf{Case 1:}  When $n=2k,$ where $k$ is a positive integer.\\
Now 
\begin{align}
&(u^{2(n+1)}-v^{2(n+1)})-(n+1)u^n v^n(u^2-v^2) \nonumber\\
&=(u^{2(2k+1)}-v^{2(2k+1)})-(2k+1)u^{2k} v^{2k}(u^2-v^2).
\end{align}
Let $a=u^2$ and $b=v^2$, so the above expression becomes
\begin{align}
&(a^{(2k+1)}-b^{(2k+1)})-(2k+1)a^{k} b^{k}(a-b)\nonumber\\
&=(a-b)\left\{\left(\sum_{i=0}^{2k} a^{2k-i}b^{i}\right)-(2k+1)a^k b^k\right\}\nonumber\\
&=(a-b)\left\{\left(\sum_{i=0}^{k-1} a^{2k-i}b^{i}\right)+\left(\sum_{i=0}^{k-1} b^{2k-i}a^{i}\right)-2ka^k b^k\right\}\nonumber\\
&=(a-b)\left\{\left(\sum_{i=0}^{k-1} a^{2k-i}b^{i}+b^{2k-i}a^{i}\right)-2ka^k b^k\right\}\nonumber\\
&=(a-b)\sum_{i=0}^{k-1} \left(a^{2k-i}b^{i}+b^{2k-i}a^{i}-2a^k b^k\right)\nonumber\\
&=(a-b)\sum_{i=0}^{k-1} a^{i} b^{i} \left(a^{2(k-i)}+b^{2(k-i)}-2a^{k-i} b^{k-i}\right)\nonumber\\
&=(a-b)\sum_{i=0}^{k-1} \left(a^{(k-i)}-b^{(k-i)}\right)^2 (ab)^i\nonumber\\
&=(u^2-v^2)\sum_{i=0}^{k-1} \left(u^{(2k-2i)}-v^{(2k-2i)}\right)^2 (uv)^{2i} \,\, {\rm as} \,\, a=u^2,~b=v^2 \nonumber\\
&=(u^2-v^2)\sum_{i=0}^{[\frac{n-1}{2}]} \left(u^{n-2i}-v^{n-2i}\right)^2 (uv)^{2i}.
\end{align}
\textbf{Case 2:}  When $n=2k-1$ where $k$ is a positive integer.\\
Now 
\begin{align}
&(u^{2(n+1)}-v^{2(n+1)})-(n+1)u^n v^n(u^2-v^2) \nonumber\\
&=(u^{4k}-v^{4k})-2ku^{2k-1} v^{2k-1}(u^2-v^2)\nonumber\\
&=(u^2-v^2)\left\{\left(\sum_{i=0}^{2k-1} u^{4k-2-2i}v^{2i}\right)-2k u^{2k-1} v^{2k-1}\right\}\nonumber\\
&=(u^2-v^2)\left\{\left(\sum_{i=0}^{k-1} u^{4k-2-2i}v^{2i}\right)+\left(\sum_{i=0}^{k-1} v^{4k-2-2i}u^{2i}\right)-2k u^{2k-1} v^{2k-1}\right\}\nonumber\\
&=(u^2-v^2)\left\{\left(\sum_{i=0}^{k-1} u^{4k-2-2i}v^{2i}+v^{4k-2-2i}u^{2i}\right)-2k u^{2k-1} v^{2k-1}\right\}\nonumber\\ 
&=(u^2-v^2)\sum_{i=0}^{k-1} \left(u^{4k-2-2i}v^{2i}+v^{4k-2-2i}u^{2i}-2 u^{2k-1} v^{2k-1}\right)\nonumber\\
&=(u^2-v^2)\sum_{i=0}^{k-1} u^{2i} v^{2i} \left(u^{4k-2-4i}+v^{4k-2-4i}-2 u^{2k-1-2i} v^{2k-1-2i}\right)\nonumber\\
&=(u^2-v^2)\sum_{i=0}^{k-1} \left(u^{2k-1-2i}-v^{2k-1-2i}\right)^2 u^{2i} v^{2i} \nonumber\\
&=(u^2-v^2)\sum_{i=0}^{[\frac{n-1}{2}]} \left(u^{n-2i}-v^{n-2i}\right)^2 u^{2i} v^{2i} .
\end{align}
\end{Proof}
\begin{Conjecture} \cite{ref3}{\label{thm2.3}}
The function $F_{n}(h):=\frac{\oint_{\Gamma_h}\,u^ny\, du}{\oint_{\Gamma_h}\,y\, du}$ is monotonically decreasing for $h\in \left(\frac{-n(n+1)^{\frac{2}{n}}}{2(n+2)},0\right).$
\end{Conjecture}
\begin{Proof}
Let
\begin{align}{\label{kp1}}
T_{n}(u):=(n+1)\frac{\int_{\delta(u)}^{u} t^n \,dt }{\int_{\delta(u)}^{u}\,dt }=u^n+u^{n-1}v+...+v^n.
\end{align}
where $v=\delta(u).$
Further  ${T_n} '(u)=\frac{1}{\Phi'(v)}\{ {f}_n(u,v)\Phi'(u)+ {f}_n(v,u)\Phi'(v)\}$, where $ {f}_n(u,v):=\sum_{k=1}^{n} kv^{k-1}u^{n-k}$.\\
Now
\begin{align}
& (u^2-v^2) \left \{ {f}_n(u,v)\Phi'(u)+ {f}_n(v,u)\Phi'(v)\right \}\nonumber\\
&=(u^2-v^2)\left\{f_n(u,w)\left(\frac{u^{n+1}}{n+1}-u\right)+f_n(u,w)\left(\frac{v^{n+1}}{n+1}-v\right)\right\} \nonumber\\
&=(u^2-v^2)\left\{\frac{1}{n+1}(u^{n+1}f_n(u,v)+v^{n+1}f_n(v,u))-(uf_n(u,v)+vf_n(v,u))\right\}\nonumber\\
&=\frac{(u^2-v^2)}{n+1}(u^{n+1}f_n(u,v)+v^{n+1}f_n(v,u))-(u^2-v^2)(uf_n(u,v)+vf_n(v,u))\nonumber\\
&=\frac{(u^2-v^2)}{n+1}(u^{n+1}f_n(u,v)+v^{n+1}f_n(v,u))-\frac{2(u^{n+2}-v^{n+2})}{(n+1)(n+2)}(uf_n(u,v)+vf_n(v,u))\nonumber\\
&~~~~~~~~~~~~~~~~~~~~~~~~~~~~~~~~~~~({\rm in \,\, view \,\, of \,\,(\ref{eqn4.2})})\nonumber\\
&=\frac{1}{(n+1)(n+2)}(n+2)(u^2-v^2)(u^{n+1}f_n(u,v)+v^{n+1}f_n(v,u))\nonumber\\
&~~~-\frac{2}{(n+1)(n+2)}(u^{n+2}-v^{n+2})(uf_n(u,v)+vf_n(v,u))\nonumber\\
&=\frac{1}{(n+1)(n+2)} (S_n(u,v)-P_n(u,v)),\label{eqn3.8}
\end{align}
where $$S_n(u,v)=(n+2)(u^2-v^2)(u^{n+1}f_n(u,v)+v^{n+1}f_n(v,u))$$ and $$P_n(u,v)=2(u^{n+2}-v^{n+2})(uf_n(u,v)+vf_n(v,u)).$$\\
Now 
\begin{align}
S_n(u,v)&= (n+2)(u^2-v^2)(u^{n+1}f_n(u,v)+v^{n+1}f_n(v,u))\nonumber\\
&=(n+2)(u^2-v^2)\left\{u^{n+1}\left(\sum_{k=1}^{n} kv^{k-1}u^{n-k}\right)+v^{n+1}\left(\sum_{k=1}^{n} ku^{k-1}v^{n-k}\right)\right\} \nonumber\\
&=(n+2)(u^2-v^2)\left\{\left(\sum_{k=1}^{n} kv^{k-1}u^{2n+1-k}\right)+\left(\sum_{k=1}^{n} ku^{k-1}v^{2n+1-k}\right)\right\} \nonumber\\
&=(n+2)\left\{\left(\sum_{k=1}^{n} kv^{k-1}u^{2n+3-k}\right)+\left(\sum_{k=1}^{n} ku^{k+1}v^{2n+1-k}\right)\right\} \nonumber\\
&~~~-(n+2)\left\{\left(\sum_{k=1}^{n} kv^{k+1}u^{2n+1-k}\right)+\left(\sum_{k=1}^{n} ku^{k-1}v^{2n+3-k}\right)\right\}  \nonumber\\
&=(n+2)\left\{\left(\sum_{k=1}^{n} kv^{k-1}u^{2n+3-k}\right)+\left(\sum_{k=1}^{n} ku^{k+1}v^{2n+1-k}\right)\right\} \nonumber\\
&~~~-(n+2)\left\{\left(\sum_{k=1}^{n} kv^{k+1}u^{2n+1-k}\right)\right\}  \nonumber\\
&~~~-(n+2)\left\{v^{2n+2} +\left(\sum_{k=2}^{n} ku^{k-1}v^{2n+3-k}\right)\right\} \nonumber\\
&=(n+2)\left\{\left(\sum_{k=1}^{n} kv^{k-1}u^{2n+3-k}\right)+\left(\sum_{k=1}^{n} ku^{k+1}v^{2n+1-k}\right)\right\} \nonumber\\
&~~~-(n+2)\left(\sum_{k=3}^{n+2} (k-2)v^{k-1}u^{2n+3-k}\right)  \nonumber\\
&~~~-(n+2)\left\{v^{2n+2} +\left(\sum_{k=0}^{n-2} (k+2)u^{k+1}v^{2n+1-k}\right)\right\} \nonumber\\
&=2(n+2)\left\{\left(\sum_{k=3}^{n} v^{k-1}u^{2n+3-k}\right)-\left(\sum_{k=1}^{n-2} u^{k+1}v^{2n+1-k}\right)\right\} \nonumber\\
&~~~+(n+2)\left(u^{2n+2}+2v u^{2n+1}-(n-1)v^n u^{n+2}-nv^{n+1} u^{n+1}\right)  \nonumber\\
&~~~+(n+2)\left((n-1)u^n v^{n+2}+n u^{n+1}v^{n+1} -v^{2n+2}-2uv^{2n+1}\right) .\label{eqn3.9}
\end{align}
Now 
\begin{align}
P_n(u,v)&= 2(u^{n+2}-v^{n+2})(uf_n(u,v)+vf_n(v,u))\nonumber\\
&=2(u^{n+2}-v^{n+2})\left\{\left(\sum_{k=1}^{n} kv^{k-1}u^{n+1-k}\right)+\left(\sum_{k=1}^{n} ku^{k-1}v^{n+1-k}\right)\right\} \nonumber\\
&=2\left\{\left(\sum_{k=1}^{n} kv^{k-1}u^{2n+3-k}\right)+\left(\sum_{k=1}^{n} ku^{n+k+1}v^{n+1-k}\right)\right\} \nonumber\\
&~~~-2\left\{\left(\sum_{k=1}^{n} kv^{n+k+1}u^{n+1-k}\right)+\left(\sum_{k=1}^{n} ku^{k-1}v^{2n+3-k}\right)\right\}\nonumber\\
&=2\left\{\left(\sum_{k=1}^{n} kv^{k-1}u^{2n+3-k}\right)+\left(\sum_{k=2}^{n+1} (n+2-k)u^{2n+3-k}v^{k-1}\right)\right\} \nonumber\\
&~~~-2\left\{\left(\sum_{k=2}^{n+1} (n+2-k)v^{2n+3-k}u^{k-1}\right)+\left(\sum_{k=1}^{n} ku^{k-1}v^{2n+3-k}\right)\right\}\nonumber\\
&=2\left\{u^{2n+2}+\left(\sum_{k=3}^{n} (n+2)v^{k-1}u^{2n+3-k}\right)+v^n u^{n+2}+(n+2)v u^{2n+1}\right\} \nonumber\\
&~~~-2\left\{u^n v^{n+2}+\left(\sum_{k=1}^{n-2} (n+2)u^{k+1}v^{2n+1-k}\right)+v^{2n+2}+(n+2)u v^{2n+1}\right\}.\label{eqn3.10}
\end{align}
Now
\begin{align}
&S_n(u,v)-P_n(u,v)\nonumber\\
&=(n+2)\left(u^{2n+2}+2v u^{2n+1}-(n-1)v^n u^{n+2}-nv^{n+1} u^{n+1}\right)  \nonumber\\
&~~~+(n+2)\left((n-1)u^n v^{n+2}+n u^{n+1}v^{n+1} -v^{2n+2}-2uv^{2n+1}\right)\nonumber\\
&~~~-2\left(u^{2n+2}+v^n u^{n+2}+(n+2)v u^{2n+1}\right)+2\left(u^n v^{n+2}+v^{2n+2}+(n+2)u v^{2n+1}\right)\nonumber\\
&~~~~~~~~~~~~~~~~~~~~~~~~~~~~~~~~~~~({\rm in \,\, view \,\, of \,\,(\ref{eqn3.9}) \,\,and\,\,(\ref{eqn3.10})})\nonumber\\
&=n(u^{2n+2}-v^{2n+2})-n(n+1)u^n v^n(u^2-v^2)\nonumber\\
&=n\left\{(u^{2n+2}-v^{2n+2})-(n+1)u^n v^n(u^2-v^2)\right\}\nonumber\\
&=n(u^2-v^2)\sum_{i=0}^{[\frac{n-1}{2}]}(u^{n-2i}-w^{n-2i})^2(uv)^{2i},\,\, {\rm in \,\, view \,\, of} \,\,Lemma\,\, \ref{kp1a}. \label{eqn3.11}
\end{align}
Using (\ref{eqn3.8}) and (\ref{eqn3.11}), we get
\begin{equation}
{f}_n(u,v)\Phi'(u)+ {f}_n(v,u)\Phi'(v)=\frac{n}{(n+1)(n+2)}\sum_{i=0}^{[\frac{n-1}{2}]}(u^{n-2i}-w^{n-2i})^2(uv)^{2i}>0 .  
\end{equation}
Now
\begin{align}
\Phi'(v)&=\frac{v^{n+1}}{n+1}-v\nonumber\\
&=v\left(\frac{v^{n}}{n+1}-1\right)<0,  \text{~as~}0<v<\sqrt[n]{n+1}.
\end{align}
This implies that 
$${T_n} '(u)=\frac{1}{\Phi'(v)}\{ {f}_n(u,v)\Phi'(u)+ {f}_n(v,u)\Phi'(v)\}<0.$$
So in view of Theorem \ref{thmc}, the function $F_{n}(h)$ is monotonically decreasing on $\left(\frac{-n(n+1)^{\frac{2}{n}}}{2(n+2)},0\right)$.
\end{Proof}

\begin{Theorem}{\label{thm2.04}}
The equation (\ref{kpr3.2}) has at most one isolated periodic traveling wave for all positive integers $n.$
\end{Theorem}
\begin{Proof}
It is shown in Section 2 that $\Gamma_h$ tends to $(\sqrt[n]{n+1},0)$ as $h \to  -\frac{n(n+1)^{\frac{2}{n}}}{2(n+2)}$ and  $\Gamma_h$ converges to the homoclinic orbit connecting $(0, 0)$ as $h\to 0.$ 
Now $A(h)=-\sqrt{c}A_0(h) \left\{\left(1+\frac{1}{c}\right) -\frac{A_n(h)}{A_0(h)} \right\}$. 
If $\frac{A_n(h)}{A_0(h)}\neq\left(1+\frac{1}{c}\right)$ for all $h\in\left(-\frac{n(n+1)^{\frac{2}{n}}}{2(n+2)},0\right)$, then $A(h)\neq0$ for all $h\in\left(-\frac{n(n+1)^{\frac{2}{n}}}{2(n+2)},0\right)$. Hence, Theorem \ref{lemma3.1} implies that the system (\ref {a01}) does not exhibit any limit cycles.

Assume that there exists a real number $c$ such that $\frac{A_n(\hat{h})}{A_0(\hat{h})}=\left(1+\frac{1}{c}\right)$  for some $\hat h\in\left(-\frac{n(n+1)^{\frac{2}{n}}}{2(n+2)},0\right)$. The monotonic nature of the ratio
$\frac{A_n(h)}{A_0(h)}$ assures that $A(h)$ has exactly one zero. Theorem \ref{lemma3.1} then suggests that the system (\ref {a01}) has a unique limit cycle, implying that the PDE (\ref{kpr3.2}) has a unique isolated periodic traveling wave solution.
\end{Proof}

\noindent \textbf{Remark:} When $n$ is even, there are two centers $(\pm\sqrt[n]{n+1},0)$ of the system (\ref{a02}). So there are two periodic annulus, one is around $(\sqrt[n]{n+1},0)$ and another is around $(-\sqrt[n]{n+1},0).$ We have proved that there is at most one limit cycle of the system (\ref {a01}) around $(\sqrt[n]{n+1},0).$ As the system (\ref{a01}) is symmetric about $u-$ axis and $y-$ axis, so there is at most one one limit cycle of the system (\ref {a01}) around $(-\sqrt[n]{n+1},0)$ by similar analysis. Clearly there is at most one isolated periodic traveling wave solution when $u(x,t)>0$ or $u(x,t)<0$.\\ \\
\textbf{Nature of wave speed:}\\ \\
Let $\Gamma_h$ be a trajectory $H(u,y)=h$ of system (\ref{a02}). \\
Now 
\begin{eqnarray}
&&H(u,y)=h\nonumber \\
&\implies& 
\frac{y^2}{2}-\frac{u^2}{2}+\frac{u^{n+2}}{(n+1)(n+2)}=h \nonumber\\
&\implies& y^2=u^2-\frac{2u^{n+2}}{(n+1)(n+2)}+2h. 
\label{eqn4.02}
\end{eqnarray}
Let $(\alpha(h),0)$ and $(\beta(h),0)$ be two intersection points of the trajectory $H(u,y)=h$ on $u-$ axis where $0\leq \alpha(h) < \beta(h).$ Also note that $H(u,y)=h$ is symmetric about $u-$ axis. 
\begin{align}
A_n(h)&=\oint_{\Gamma_ h}\,u^n y\, du\nonumber\\
&=2\int_{\alpha(h)}^{\beta(h)} \,u^n y\, du \nonumber\\
&:=2 J_n(h).
\end{align}

\begin{Lemma} {\cite{ref1}} \label{lem2}
Let $B(p,q)=\int_{0}^{1} x^{p-1} (1-x)^{q-1} dx ,$ $p>0,$ $q>0$ be the Beta function. Then 
\begin{equation*}
J_0(0)= \frac{1}{n} \left(\frac{(n+1)(n+2)}{2}\right)^{\frac{2}{n}} B\left(\frac{3}{2},\frac{2}{n}\right),
\end{equation*}

\begin{equation*}
J_0(0)= \frac{1}{n} \left(\frac{(n+1)(n+2)}{2}\right)^{\frac{n+2}{n}} B\left(\frac{3}{2},\frac{n+2}{n}\right),    
\end{equation*}
and 
\begin{equation*}
 \frac{J_n(0)}{J_0(0)}=\frac{2(n+1)(n+2)}{3n+4} ,\, \, \,\lim_{h \to -\frac{n(n+1)^{\frac{2}{n}}}{2(n+2)}} \frac{J_n(h)}{J_0(h)}= n+1. 
\end{equation*}
\end{Lemma}
Now Lemma \ref{lem2} indicates that $\frac{A_n(0)}{A_0(0)}=\frac{2J_n(0)}{2J_0(0)}=\frac{2(n+1)(n+2)}{3n+4}>1,$ and\\
\begin{equation*}
\lim_{h \to -\frac{n(n+1)^{\frac{2}{n}}}{2(n+2)}} \frac{A_n(h)}{A_0(h)}= \lim_{h \to -\frac{n(n+1)^{\frac{2}{n}}}{2(n+2)}} \frac{2J_n(h)}{2J_0(h)}=n+1 .   
\end{equation*}
As $F_{n}(h)=\frac{A_n(h)}{A_0(h)}$ is monotonically decreasing on the interval $\left(\frac{-n(n+1)^{\frac{2}{n}}}{2(n+2)},0\right),$ we have
\begin{equation} \label{eqn3.16}
   \frac{2(n+1)(n+2)}{3n+4}<\frac{A_n(h)}{A_0(h)}<n+1. 
\end{equation}
Next we quote a part of result of Theorem 1 of \cite{ref1} to discuss about limit wave speed of periodic traveling waves.
\begin{Lemma} \cite{ref1}
There exists $\hat \epsilon>0,$ then for $\epsilon \in (0,\hat \epsilon)$ and $h \in \left(\frac{-n(n+1)^{\frac{2}{n}}}{2(n+2)},0\right),$ wave speed $c=c(\epsilon,h)$ of periodic traveling wave is a smooth function of $h$ and $\epsilon.$ Also $c(\epsilon,h)$ converges to $c_{0}(h)$ as $\epsilon \to 0,$ where $c_0(h)$ is a smooth function for $h \in \left(\frac{-n(n+1)^{\frac{2}{n}}}{2(n+2)},0\right).$
\end{Lemma}
\noindent The function $c_{0}(h)$ is called limit wave speed. We have below Theorem for an important property of limit wave speed $c_{0}(h).$
\begin{Theorem}{\label{thm2.4}}
 Limit wave speed $c_0(h)$ is monotonically increasing on $\left(\frac{-n(n+1)^{\frac{2}{n}}}{2(n+2)},0\right).$
 \end{Theorem}

\begin{Proof}
From equation (\ref{eqn2.11}), we know that when $\epsilon \to 0,$ there exists a periodic wave if 
\begin{eqnarray}
&&\frac{A_n(h)}{A_0(h)}=1+\frac{1}{c_{0}(h)} \nonumber \\
&\implies& 
c_{0}(h)=\frac{1}{\frac{A_n(h)}{A_0(h)}-1} \label{eqn3.17}
\end{eqnarray}
In view of proof of Conjecture \ref{thm2.3}, $\frac{A_n(h)}{A_0(h)}$ is monotonically decreasing on $\left(\frac{-n(n+1)^{\frac{2}{n}}}{2(n+2)},0\right).$ This indicates that $c_{0}(h)$ is monotonically increasing on $\left(\frac{-n(n+1)^{\frac{2}{n}}}{2(n+2)},0\right)$ in view of (\ref{eqn3.17}). Moreover,
\begin{eqnarray}
&&\frac{2(n+1)(n+2)}{3n+4}<\frac{A_n(h)}{A_0(h)}<n+1 \text{~~for~} h \in \left(\frac{-n(n+1)^{\frac{2}{n}}}{2(n+2)},0\right)\nonumber \\
&&~~~~~~~~~({\rm in \,\, view \,\, of \,\,(\ref{eqn3.16})})\nonumber\\
&\implies& 
\frac{2(n+1)(n+2)}{3n+4}<1+\frac{1}{c_{0}(h)}<n+1 \nonumber\\
&\implies& \frac{2n^2+3n}{3n+4}<\frac{1}{c_{0}(h)}<n  \nonumber\\
&\implies& \frac{1}{n}<c_{0}(h)<\frac{3n+4}{2n^2+3n}.
\label{eqn4.002}
\end{eqnarray}
Also 
\begin{equation}
  \lim_{h \to -\frac{n(n+1)^{\frac{2}{n}}}{2(n+2)}} c_{0}(h)=\frac{1}{n}, \, \, \text{and} \, \, c_{0}(0)=\frac{3n+4}{2n^2+3n}.
\end{equation}
\end{Proof}
\section{Numerical simulation}
In this section, we present numerical simulations to verify the theoretical results discussed earlier. 

Let us consider $n=5,$ $\epsilon=0.1$ and the initial values $(u,y)=(\frac{1}{2},0)$. Clearly $H(\frac{1}{2},0)=-\frac{671}{5376}:=h^*$ and the corresponding trajectory of the unperturbed system (\ref{a02}) is $\Gamma_{h^*}: H(u,y)=\frac{y^2}{2}-\frac{u^2}{2}+\frac{u^7}{42}=-\frac{671}{5376}$. 
By direct computation, we get  $\frac{A_5(h^*)}{A_0(h^*)}=\frac{\oint_{\Gamma_{h^*}}\,u^5y\, du}{\oint_{\Gamma_{h^*}}\,y\, du}=4.88851$. Let us choose $c=\frac{100000}{388851}$. This clearly gives $A(h^*)=0.$ \\ 

\begin{figure}[htbp] 
  \centering
      \includegraphics[width=1.08\linewidth]{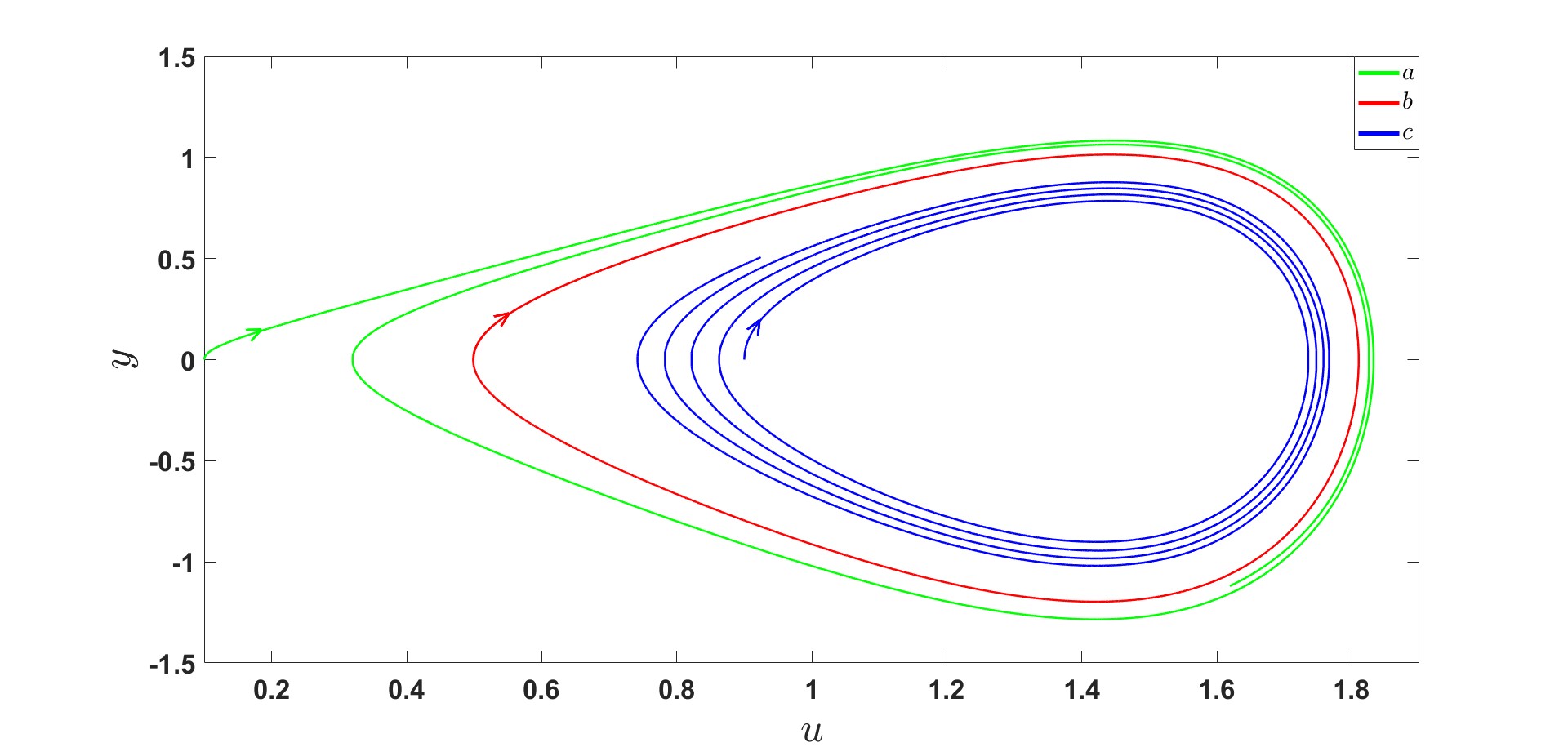}
  \caption{Trajectory of system (\ref{a01}) for $\epsilon=0.1,$ $n=5,$ $c=\frac{100000}{388851}$ with initial condition $(a)$ $u=0.1$ and $y=0,$ $(b)$ $u=0.49885$ and $y=0,$ $(c)$ $u=0.9$ and $y=0.$}
  \label{kpr005}
\end{figure}

Observations from numerical simulations reveal that:

\begin{itemize}
    \item [\rm(i)]The perturbed system (\ref{a01}) exhibits an inward spiraling trajectory starting from the initial point $(0.1, 0)$ (as shown in Figure \ref{kpr005}-$a$), while the trajectory with an initial point of $(0.9, 0)$ spirals outward (as shown in Figure \ref{kpr005}-$c$). This indicates the presence of a stable limit cycle.
    
\item [\rm(ii)] The trajectory of the perturbed system (\ref{a01}) beginning from the initial condition $(0.49885, 0)$ closely approximates a closed trajectory (see Figure \ref{kpr005}-$b$).
\end{itemize} 

\begin{figure}[htbp] 
  \centering
      \includegraphics[width=1.08\linewidth]{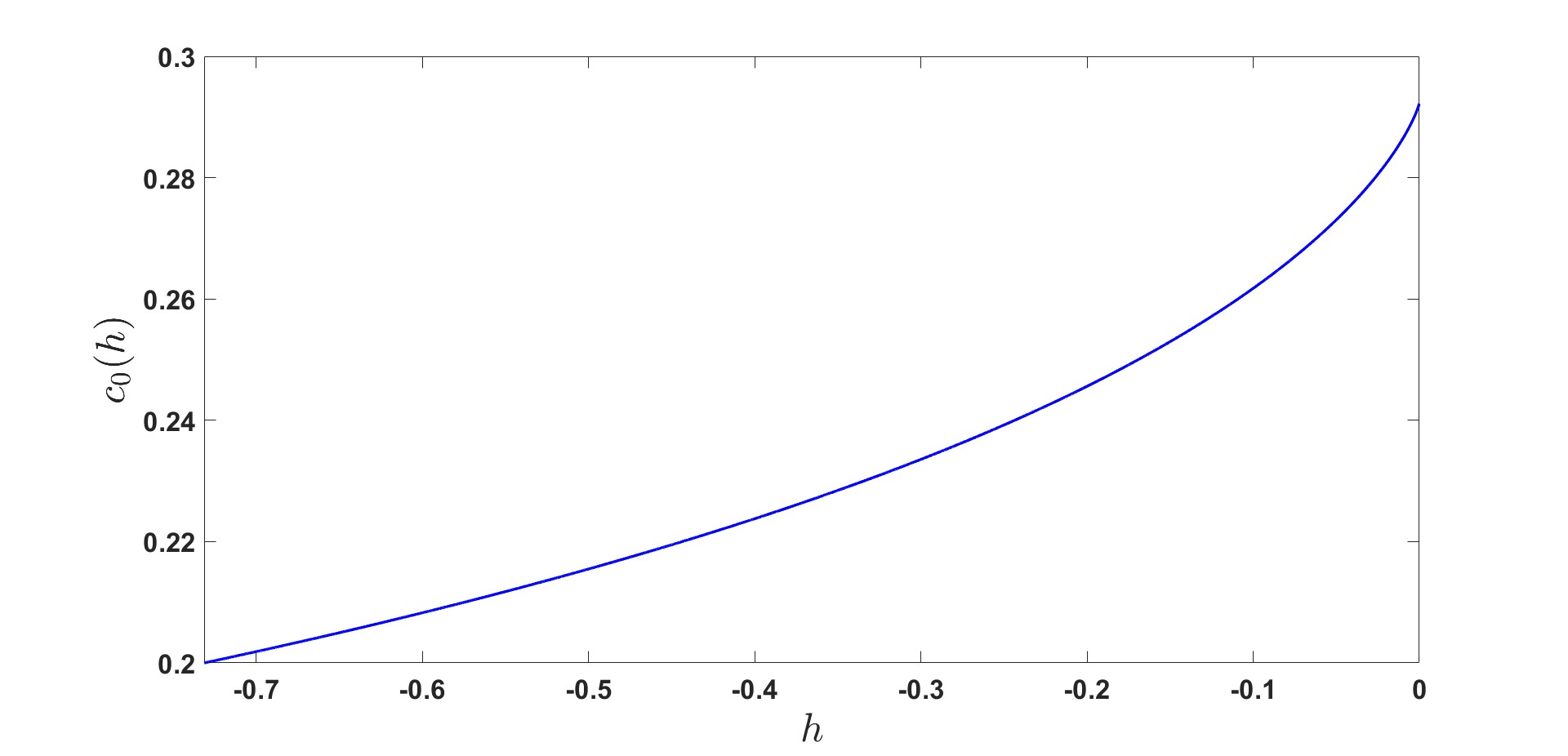}
  \caption{Graphical representation of the limit wave speed $c_{0}(h)$ for $n=5.$}
  \label{kpr004}
\end{figure}

When $n=5,$ and $\epsilon \to 0,$ the theoretical results indicate that $c_{0}(h)$ is increasing over the interval $\left(-\frac{5\sqrt[5]{36}}{14},0\right),$ which aligns with the numerically plotted curve shown in Figure \ref{kpr004}.

\section{Conclusion}
Our analysis in this paper centers on the monotonic behavior of the limit wave speed for periodic wave solutions of the perturbed generalized KdV equation. The methodologies introduced include the application of geometric singular perturbation theory for reducing singularly perturbed systems to regular perturbed systems. The conjecture proposed by Ouyang et al.\cite{ref3} has been proven in this work. In addition, we showed that the limit wave speed $c_{0}(h)$ of the periodic wave is monotone for $h \in \left(\frac{-n(n+1)^{\frac{2}{n}}}{2(n+2)},0\right),$ addressing a previously open question. Furthermore, we proved that there is at most one isolated periodic wave solution when when $u(x,t)>0$ or $u(x,t)<0.$ We establish the criteria for selecting the wave speed of periodic wave solution for equation (\ref{kpr3.2}). Numerically, for $n=5,$ we have shown that a stable limit cycle exists for the system (\ref{a01}), and we plotted the graph of $c_{0}(h)$ to validate the theoretical results.\\ \\
\noindent\textbf{Author Contributions:} The manuscript was written and reviewed by both authors.\\
\textbf{Data Availability:} This document does not have any associated manuscript.
\section*{Declarations}
\textbf{Conflict of interest:} The authors state that this research was conducted without any conflicts of interest.

\end{document}